\documentclass[letterpaper, 10 pt, conference]{ieeeconf}
\IEEEoverridecommandlockouts                              

\usepackage{graphicx} 
\usepackage{url}
\usepackage{ifthen}      

\usepackage{mathtools} 
\newboolean{techrep}
\setboolean{techrep}{false}
\setboolean{techrep}{true}

\makeatletter\let\NAT@parse\undefined\makeatother 
\usepackage[numbers,sort&compress]{natbib} 

\usepackage{url}
\usepackage[draft,ieee,fancybb]{jphmacros2e}

\ifthenelse{\boolean{techrep}}{
  \newcommand{\revision}{}
}{
  \newcommand{\revision}{\color{blue}}
}

\title{\LARGE \bf Markov Chain
  Monte Carlo 
  for Koopman-based Optimal Control\ifthenelse{\boolean{techrep}}{:
    Technical Report}{}}

\author{Jo\~ao Hespanha$^{1}$ \and
  Kerem \c{C}amsar\i$^{2}$
  \thanks{*This material is based upon work supported by the
    U.S. Office of Naval Research MURI grant No.~N00014-23-1-2708 and
    by the National Science Foundation grant
    No.~2229876.}
  \thanks{Univ.~of California, Santa Barbara, USA;
    $^{1}$\protect\url{hespanha@ece.ucsb.edu},
  $^{2}$\protect\url{camsari@ece.ucsb.edu}}%
}

\begin{document}

\maketitle
\thispagestyle{empty}
\pagestyle{empty}

\begin{abstract}
  We propose a Markov Chain Monte Carlo (MCMC) algorithm
  based on Gibbs sampling with parallel tempering to solve nonlinear
  optimal control problems. The algorithm is applicable to nonlinear
  systems with dynamics that can be approximately represented by a
  finite dimensional Koopman model, potentially with high dimension.
  {\revision This algorithm exploits linearity of the Koopman
    representation to achieve significant computational saving for
    large lifted states. We use a video-game to illustrate the use of
    the method.}
\end{abstract}


\section{INTRODUCTION}


\ifthenelse{\boolean{techrep}}{ This paper address the
    optimal control of nonlinear systems that have reasonably accurate
    finite-dimensional representations in terms of their Koopman
    operator~\cite{Koopman1931}.}{} {\revision While \citeauthor{Koopman1931}'s
  pioneering work is almost a century old, its use as a practical tool
  to model complex dynamics is much more recent and only became
  practical when computational tools became available for the analysis
  of systems with hundreds to thousands of dimensions. The use of Koopman
  models for control is even more recent, but has attracted
  significant attention in the last few years
  \cite{BruntonBruntonProctorKutz2016,KordaMezic2018,ProctorBruntonKutz2018,MauroySusukiMezic2020,BevandaSosnowskiHirche2021,OttoRowley2021}\nocite{Mezic2005}.}

\smallskip

The linear structure of the Koopman representation
permits very efficient solutions for optimal control when the lifted
dynamics are linear on the control input, as
in~\cite{BruntonBruntonProctorKutz2016,KordaMezic2018,ProctorBruntonKutz2018}. However,
linearity in the control severely limits the class of applicable
dynamics. Bilinear representations are more widely
applicable~\cite{MauroySusukiMezic2020,BevandaSosnowskiHirche2021,OttoRowley2021},
but are also harder to control.

\smallskip

When the set of admissible control inputs is finite, the Koopman
representation can be viewed as a switched linear system, where the
optimal control 
selects, at each time step, one out of several admissible dynamics
\cite{BlischkeHespanhaDec2023}. The optimal control of switched system
it typically difficult
\cite{bengea2005optimal,xu2004optimal}, 
but when the optimization criterion is linear in the lifted state, the
dynamic programming cost-to-go is concave and piecewise linear, with a
simple representation in terms of a minimum over a finite set of
linear functions. This observation enabled the design of efficient
algorithms that combine dynamic programming with dynamic
pruning~\cite{BlischkeHespanhaDec2023}.

\smallskip

{\revision This paper exploits the structure of the Koopman
  representation to develop efficient Markov Chain Monte Carlo (MCMC)
  sampling method for optimal control. Following the pioneering work
  of
  \cite{Kirkpatrick1983optimization,Cerny1985thermodynamical,GemanGeman1984stochastic},
  we draw samples from a Boltzmann distribution with energy
  proportional to the criterion to minimize.
%
%
The original use of MCMC methods for combinatorial optimization relied
on a gradual decrease in temperature, now commonly known as
\emph{simulated annealing}, to prevent the chain from getting trapped
into states that do not minimize energy. An alternative approach
relies on the use of multiple replicas of a Markov chain, each
generating samples for a different temperature. The introduction of
multiple replicas of a Markov chain to improve the mixing time can be
traced back to \cite{Swendsen1986replica}. The more recent form of
\emph{parallel tempering} (also known as \emph{Metropolis–coupled
  MCMC}, or \emph{exchange Monte Carlo}) is due to
\cite{Geyer1991markov,Hukushima1996exchange}.}

\smallskip

Our key contribution is an MCMC algorithm that combines Gibbs sampling
\cite{GemanGeman1984stochastic,Gelfand1990sampling} with parallel
tempering to solve the switched linear optimizations that arise from
the Koopman representation of nonlinear optimal control problems. This
algorithm is computationally efficient due to the combination of two
factors: (i) the linear structure of the cost function enables the
full variable sweep need for Gibbs sampling to be performed with
computation that scales linearly with the horizon length $T$ and (ii)
parallel tempering can be fully parallelized across computation cores,
as noted in~\cite{Earl2005parallel}. {\revision With regard to (i),
  the computational complexity of evaluation Gibbs' conditional
  distribution for each optimization variable is of order
  $T n_\psi^2$, where $n_\psi$ denotes the size of the lifted state. Our
  algorithm computes the conditional distributions \emph{for all $T$
    variables} in a full Gibbs' sweep with computation still just of
  order $T n_\psi^2$.} With regard to (ii), while here we only explore
parallelization across CPU cores, {\revision in the last few years}
hardware parallelization using GPUs and FPGAs has achieved orders of
magnitude increase in the number of samples generated with MCMC
sampling~\cite{mohseni2022ising,aadit2022massively}.

\smallskip

The remaining of this paper is organized as follows:
Section~\ref{se:Koopman} shows how a nonlinear control problem can be
converted into a switching linear systems optimization, using the
Koopman operator. Section~\ref{se:MCMC} provides basic background on
MCMC, Gibbs sampling, and parallel tempering. Our optimization
algorithm is described in Section~\ref{se:MCMC-Koopman} and its use is
illustrated in Section~\ref{se:numerical} in the context of a video
game. \ifthenelse{\boolean{techrep}}{}{While this paper is
  self-contained, details of some of the algebraic derivations are
  omitted, but can be found in~\cite{HespanhaMar2024b_arXiv}.}

\section{OPTIMAL CONTROL OF KOOPMAN MODELS}\label{se:Koopman}

Given a discrete-time nonlinear system of the form
\begin{align}\label{eq:nonlinear-process}
  x_{t+1}=f(x_t,u_t), \quad
  \forall t\in\scr{T},\;
  x_t\in\scr{X}, \;
  u_t\in\scr{U}_t,
\end{align}
with the time $t$ taking values over $\scr{T}\eqdef\{1,\dots,T\}$, our
goal is to solve a final-state optimal control problem of the form
\begin{align}\label{eq:nonlinear-minimization}
  J^*\eqdef \min_{u\in \scr{U}} J(u), \quad
  J(u)\eqdef g(x_T),
\end{align}
where
$u\eqdef(u_1,\dots,u_T)\in \scr{U}\eqdef\scr{U}_1\times\cdots\times \scr{U}_T$ denotes
the control sequence to be optimized, which we assume finite but
potentially with a large number of elements.

\smallskip

For each input $u\in\scr{U}_t$, $t\in\scr{T}$, the Koopman operator
$K_u$ for the system~\eqref{eq:nonlinear-process} operates on the
linear space of functions $\scr{F}$ from $\scr{X}$ to $\R$ and
is defined by
\begin{align*}
  \varphi(\cdot)\in\scr{F}  \mapsto \varphi(f(\cdot,u))\in \scr{F}.
\end{align*}
Assuming there is a finite dimensional linear subspace
$\scr{F}_\mrm{inv}$ of $\scr{F}$ that is invariant for every Koopman
operator in the family
$\{K_u: \forall u\in\scr{U}_t, \; t\in\scr{T}\}$, there is an associated family
of matrices $\{A(u)\in\R^{n_\psi\times n_\psi}: \forall u\in\scr{U}_t, \; t\in\scr{T}\}$ such that
\begin{align}\label{eq:Koopman-model}
  \psi_{t+1}= A(u_t) \psi_t, \quad
  \forall t\in\scr{T},\;
  \psi_t\in\R^{n_\psi}, \;
  u_t\in\scr{U}_t,
\end{align}
where
$\psi_t \eqdef \matt{\varphi_1(x_t) & \cdots & \varphi_{n_\psi}(x_t)}' \in\R^{n_\psi}$ and the
functions $\{\varphi_1(\cdot),\dots,\varphi_{n_\psi}(\cdot)\}$ form a basis for
$\scr{F}_\mrm{inv}$
\cite{FolkestadBurdick2021,BlischkeHespanhaDec2023,HaseliCortes2023}. If in addition, the
function $g(\cdot)$ in~\eqref{eq:nonlinear-minimization} also belongs to
$\scr{F}_\mrm{inv}$, there also exists a row vector $c\in\R^{1\times n_\psi}$ such that
\begin{align}\label{eq:Koopman-cost}
  g(x_T) = c \psi_T,
\end{align}
which enables re-writing the optimization
criterion~\eqref{eq:nonlinear-minimization} as
\begin{align}\label{eq:Koopman-criterion}
  J(u)\eqdef c \psi_T = c A(u_T)\cdots A(u_1) x_1.
\end{align}
We can thus view the original optimal control problem for the
nonlinear system~\eqref{eq:nonlinear-process} as a switched linear
control problem~\cite{BlischkeHespanhaDec2023}. {\revision Note that
  it is always possible to make sure that \eqref{eq:Koopman-cost} hold
  for some vector $c$ by including $g(x_t)$ as one of the entries of
  $\psi_t$.}

\smallskip

It is generally not possible to find a finite-dimensional subspace
$\scr{F}_\mrm{inv}$ that contains $g(\cdot)$ and is invariant for every
$\{K_u:u\in\scr{U}_t,\,t\in\scr{T}\}$. {\revision Instead, we typically
  work with a finite dimensional space for
  which~\eqref{eq:Koopman-model}--\eqref{eq:Koopman-cost} hold up to
  some error. However, to make this error small, we typically need to
  work with high-dimensional subspaces, i.e., large values of
  $n_\psi$. \revision Quantifying the impact on the error of the lifted
  space dimension $n_\psi$ and of the size of the dataset used to learn
  the Koopman dynamics is challenging, but notable progress has been
  reported
  in~\cite{HaseliCortes2023,NuskePeitzPhilippSchallerWorthmann2023}. Much
  work remains, e.g., in quantifying errors for systems with isolated
  limit sets~\cite{bakker2019learning,liu2023properties}.}

\section{MCMC METHODS FOR OPTIMIZATION}\label{se:MCMC}

To solve a combinatorial minimization of the form
\begin{align}\label{eq:optimization-discrete}
  J^*\eqdef \min_{u\in\scr{U}} J(u).
\end{align}
over a finite set $\scr{U}$, it is convenient to consider the following
\emph{Boltzmann distribution} with energies $J(u)$, $u\in\scr{U}$:
\begin{align}\label{eq:Boltzmann}
  p(u;\beta) &\eqdef \frac{e^{-\beta J(u)}}{Q(\beta)}, \; \forall u\in\scr{U},&
  & Q(\beta) \eqdef \sum_{\bar u\in\scr{U}} e^{-\beta J(\bar u)},
\end{align}
for some constant $\beta\ge 0$. For consistency with statistical mechanics,
we say that values of $\beta$ close to zero correspond to high
temperatures, whereas large values correspond to low
temperatures. \ifthenelse{\boolean{techrep}}{ The normalization
  function $Q(\beta)$ is called the \emph{canonical partition function}.}

\smallskip

For $\beta=0$ (high temperature), the Boltzmann distribution becomes
uniform and all $u\in\scr{U}$ are equally probable. However, as we
increase $\beta$ (lower the temperature) all probability mass is
concentrated on the subset of $\scr{U}$ that minimizes the energy
$J(u)$.  \ifthenelse{\boolean{techrep}}{ This can be seen by noting
  that the ratio of the probability between two $u,\bar u\in\scr{U}$ is
  given by
\begin{align*}
  \frac{p(u;\beta)}{p(\bar u;\beta)} = e^{-\beta(J(u)-J(\bar u))}
  \xrightarrow[\beta\to+\infty]{}
  \begin{cases}
    0 & J(u)>J(\bar u)\\
    1 & J(u)=J(\bar u)\\
    \infty& J(u)<J(\bar u),
  \end{cases}
\end{align*}
which shows that the value with higher energy (and higher cost) will
eventually never be selected.

\smallskip }{} This motivates a procedure to solve
\eqref{eq:optimization-discrete}: draw samples from a random variable
$\mbf u_\beta$ with Boltzmann distribution given by~\eqref{eq:Boltzmann}
with $\beta$ sufficiently high (temperature sufficiently low) so that all
samples correspond to states with the minimum energy/cost.

\ifthenelse{\boolean{techrep}}{
\smallskip

The following result confirms that the probability $\delta$ that we can
extract $N$ samples from the Boltzmann distribution with energies much
larger than the minimum $J^*$ can be made arbitrarily small by
increasing $N$ in proportion to just the \emph{logarithm} of $\delta$,
provided that $\beta$ is sufficiently high (temperature is sufficiently
low).
\begin{lemma}\label{le:Boltmann-sample-complexity}
  Let $\{\mbf u[1],\dots,\mbf u[N]\}$ be independent samples from the
  Boltzmann distribution~\eqref{eq:Boltzmann}. For every
  $\epsilon\ge 0$, $\delta>0$
  \begin{align*}
    \Prob\big(J(\mbf u[i])\ge J^*+\epsilon,\;\forall i\big) \le\delta,
  \end{align*}
  provided that
  \begin{align*}
    \beta &\ge \frac{\log(|\scr{U}|/|\scr{U}^*|)}{\epsilon}, &
    N &\ge \frac{\log(1/\delta)}{\beta \epsilon -\log(|\scr{U}|/|\scr{U}^*|)},
  \end{align*}
  where
  $\scr{U}^*=\argmin_{u\in\scr{U}} J(u)\eqdef \{u\in\scr{U}:J(u)= J^*\}$.
  \frqed
\end{lemma}
The number of samples $N$ required to achieve a small probability
$\delta$ is a function of the fraction $|\scr{U}^*|/|\scr{U}|$ of elements
of $\scr{U}$ that minimize $J(u)$. The result is particularly
interesting in that the magnitude of $\beta$ and number of samples $N$
depends on the \emph{logarithm} of this fraction, meaning that the
size of $\scr{U}$ could grow exponentially with respect to some
scaling variable (like the time horizon in optimal control) and yet
$\beta$ and $N$ would only need to grow linearly. However, the optimism
arising from this bound needs to be tempered by the observation that
this result assumes that the samples are independent and this is
typically hard to achieve for large values of $\beta$, as discussed below.
}{}

\subsection{Markov Chain Monte Carlo Sampling}\label{se:MCMC-sampling}

Consider a \emph{discrete-time 
  Markov chain} $\{\mbf u[1],\mbf u[2],\dots\}$ on a finite set
$\scr{U}$, with \emph{transition
  probability}
\begin{align*}
  p(u' | u)=\Prob\big(\mbf u[k+1]=u'\,|\,\mbf u[k]=u),
  \;\forall u',u\in\scr{U},\, k\ge1.
\end{align*}
Combining the probabilities of all possible {\revision realizations}
of $\mbf u[k]$ in a row vector \ifthenelse{\boolean{techrep}}{
\begin{align*}
  p[k]\eqdef \big[ \Prob(\mbf u[k] =u) \big]_{u\in\scr{U}}
\end{align*}
}{ $p[k]$ } and organizing the values of the transition probabilities
$p(u' |u)$ as a \emph{transition
  matrix}\ifthenelse{\boolean{techrep}}{
\begin{align*}
  P \eqdef \Big[ p(u' | u) \Big]_{u,u'  \in\scr{U}},
\end{align*}
}{ $P$,} with one $u'\in\scr{U}$ per column and one $u\in\scr{U}$ per row,
{\revision enables} us to express the evolution of the $p[k]$
as
\begin{align*}
  p[k+K] = p[k] P^K, \quad \forall k,K\ge 1.
\end{align*}
A Markov chain is called \emph{regular} if that there exists an
integer $N$ such that all entries of $P^N$ are strictly positive. In
essence, this means that any $\bar u\in \scr{U}$ can be reached in $N$
steps from any other $\tilde u\in \scr{U}$ through a sequence of
transitions with positive probability $p(u' |u)>0$,
$u,u'\in\scr{U}$. {\revision The following result adapted
  from~\cite[Chapter IV]{KemenySnell1976} is the key property behind
  MCMC sampling:}
\begin{theorem}[Fundamental Theorem of Markov Chains]
  \label{th:fundamental-Markov-chains}
  For every regular Markov chain on a finite set $\scr{U}$ and
  transition matrix $P$, there exists a vector $\pi$ such that:
  \begin{enumerate}
  \item {\revision The Markov chain is \emph{geometrically ergodic},
      meaning that there exists constants $c> 0,\lambda\in[0,1)$ such that}
    \begin{align}\label{eq:exponential}
      \|p[k] P^K-\pi\| \le c \lambda^K, \quad \forall k,K\ge 0.
    \end{align}

  \item The vector $\pi$ is called the \emph{invariant distribution} and
    is the unique solution to the \emph{global balance equation}:
    \begin{align}\label{eq:global-balance-P}
      &\pi P = \pi, &
      &\pi \sbf 1=1, &
      & \pi \ge 0.
    \end{align}
  \end{enumerate}
\end{theorem}
\smallskip

{\revision To use an MCMC method to draw samples from a desired distribution
$p(u;\beta)$, we construct a discrete-time Markov Chain that is regular
and with an invariant distribution $\pi$ that matches $p(u;\beta)$. We
then ``solve'' our sampling problem by
only accepting a subsequence of samples
$\mbf u[K+1], \mbf u[2K+1], \dots$ with $K$ ``sufficiently large'' so
that the distribution $p[K+1]=p[1] P^K$ of the sample $\mbf u[K+1]$
satisfies
\begin{align}\label{eq:K-bound}
  \|p[1] P^K-\pi\|
  \le \epsilon  &&
  \xLeftarrow{\text{\eqref{eq:exponential}}}
  &&  K\ge \frac{\log c+\log\epsilon^{-1}}{\log\lambda^{-1}}.
\end{align}
for some ``sufficiently small'' $\epsilon$. Such $K$ guarantees that the
samples $\mbf u[K+1], \mbf u[2K+1], \dots$ may not quite have the
desired distribution, but are away from it by no more than
$\epsilon$.}  \ifthenelse{\boolean{techrep}}{
Two samples separated by $K$ are not quite independent, but
``almost''. Independence between $\mbf u[K+1]$ and $\mbf u[2K+1]$ would
mean that the following two conditional distributions must have the
same value
\begin{align*}
  &\Prob(\mbf u[2K+1]\,|\,\mbf u[K+1]=u) = p_u P^K, \\
  &\Prob(\mbf u[2K+1]\,|\,\mbf u[K+1]=\bar u)= p_{\bar u} P^K,
\end{align*}
for every two distributions $p_u$ and $p_{\bar u}$ that place all the
probability weight at $u$ and $\bar u\in\scr{U}$, respectively. In
general these conditional probabilities will not be exactly equal, but
we can upper bound their difference by
\begin{align*}
  \|p_u P^K-p_{\bar u} P^K\|
  \le \|p_u P^K-\pi\|+\|\pi-p_{\bar u} P^K\| \le 2\epsilon.
\end{align*}
While we do not quite have independence, this shows that any two
conditional distributions of $\mbf u[K+1], \mbf u[2K+1], \dots$ are within
$2\epsilon$ of each other.
\subsubsection{Balance}
}{ {\revision Since $\epsilon^{-1}$ appears in~\eqref{eq:K-bound} ``inside''
    a logarithm, we can get $\epsilon$ extremely small without having to
    increase $K$ very much.} The factor $1/\log\lambda^{-1}$
  in~\eqref{eq:K-bound} is often called the \emph{mixing time} of the
  Markov chain and is a quantity that should be small to minimize the
  number of ``wasted samples.''

  \smallskip
}

To make sure that the chain's invariant distribution $\pi$ matches a
desired sampling distribution $p(u;\beta)$, $u\in\scr{U}$, we need the
latter to satisfy the \emph{global balance
  equation}~\eqref{eq:global-balance-P}\ifthenelse{\boolean{techrep}}{,
  which can be re-written in non-matrix form as
\begin{align}\label{eq:global-balance}
  p(u';\beta) = \sum_{u \in\scr{U}} p(u'|u)p(u;\beta), \quad\forall u'\in\scr{U}.
\end{align}
}{.}
{\revision A sufficient (but not necessary) condition for global balance is
\emph{detailed balance}~\cite{Brooks1998},} which instead
  asks for
\begin{align}\label{eq:detailed-balance}
  p(u|u')p(u';\beta)=p(u'|u)p(u;\beta), \quad\forall u',u \in\scr{U}.
\end{align}

\ifthenelse{\boolean{techrep}}{
\subsubsection{Mixing Times}

We saw above that to obtain (approximately) independent samples with
the desired distribution we can only use one out of $K$ samples from
the Markov chain with $K$ satisfying~\eqref{eq:K-bound}.

\smallskip

Since $\epsilon^{-1}$ appears in~\eqref{eq:K-bound} ``inside'' a logarithm,
we can get $\epsilon$ extremely small without having to increase $K$ very
much. However, the dependence of $K$ on $\lambda$ can be more problematic
because $\lambda\in[0,1)$ can be very close to 1. The multiplicative factor
$\frac{1}{\log\lambda^{-1}}$ is typically called the \emph{mixing time} of
the Markov chain and is a quantity that we would like to be small to
minimize the number of ``wasted samples.''

\smallskip

Computing the mixing time is typically difficult, but many bounds for
it are available~\cite{RothblumTan1985}.  One of the simplest bounds
due to Hoffman states that
\begin{multline}\label{eq:Hoffman}
  \lambda \le \frac{p_{\max}-p_{\min}}{p_{\max}+p_{\min}}
  \\
  \imply
  \frac{1}{\log \lambda^{-1}}
  \le  \frac{1}{1-\lambda}
  \le \frac{p_{\max}+p_{\min}}{2p_{\min}}
\end{multline}
where
\begin{align*}
  &p_{\max}\eqdef 
  \sum_{u' \in\scr{U}} \max_{u\in\scr{U}} p(u'|u), &
  &p_{\min}\eqdef 
  \sum_{u' \in\scr{U}} \min_{u\in\scr{U}} p(u'|u).
\end{align*}
}{}

\subsection{Gibbs Sampling}

\ifthenelse{\boolean{techrep}}{
\emph{Gibbs sampling} is an MCMC sampling algorithm that generates a
Markov chain with a desired multivariable distribution. This method
assumes that, while sampling from the joint distribution is difficult,
it is easy to sample from conditional distributions.

\smallskip

We start }{ \emph{Gibbs sampling} starts }%
with a function $p(u;\beta)$, $u\in\scr{U}$ that defines a desired
multi-variable joint distribution up to a normalization constant. We
use the subscript notation $u_t$ to refer to the variable
$t\in\scr{T}$ in $u$ and $u_{-t}$ to refer to the remaining
variables. The algorithm operates as follows:

\smallskip

\begin{algorithm}[Gibbs sampling]~\small
  \label{alg:Gibbs}
\begin{senumerate}{0.5ex}
\item Pick arbitrary $\mbf u[1]\in\scr{U}$.
\item\label{en:loop} Set $k=1$ and repeat until enough samples are collected:
  \begin{itemize}
  \item \emph{Variable sweep}: For each $t\in\scr{T}$:

    \begin{itemize}
    \item Sample $\mbf u_t[k+1]$ with the desired conditional
      distribution of $\mbf u_t$, given $\mbf u_{-t}[k]$:
      \begin{align}\label{eq:Gibbs-i}
        \frac{p(u_t,\mbf u_{-t}[k];\beta)}{\sum_{\bar u_t} p(\bar u_t,\mbf u_{-t}[k];\beta)},
      \end{align}

    \item Set $\mbf u_{-t}[k+1]=\mbf u_{-t}[k]$ and increment $k$.\frqed
    \end{itemize}
  \end{itemize}
\end{senumerate}
\end{algorithm}

\smallskip

\ifthenelse{\boolean{techrep}}{
\subsubsection{Balance}

It is straightforward to show that the Gibbs transition probability
for a sample corresponding to the update of each variable
$\mbf u_t[k]$, $t\in\scr{T}$ satisfies the detailed balance
equation~\eqref{eq:detailed-balance} for the desired distribution
$p(u;\beta)$:
\begin{align*}
  p_t(u'|u)p(u;\beta)
  &=p_t(u|u')p(u';\beta), \quad\forall u,u' \in\scr{U}
\end{align*}
which means that we also have global balance.  Denoting by $P_t$ the
corresponding transition matrix, this means that
\begin{align*}
  \pi_\beta P_t=\pi_\beta, \quad\forall t\in\scr{T},
\end{align*}
where $\pi_\beta$ denotes the vector of probabilities associated with the
desired distribution $p(u;\beta)$. A full variable sweep has transition
matrix given by
\begin{align*}
  P_\mrm{sweep}\eqdef P_1 P_2\cdots P_T,
\end{align*}
which also satisfies global balance because
\begin{align*}
  \pi_\beta P_\mrm{sweep}
  =\pi_\beta P_1 P_2\cdots P_T
  =\pi_\beta P_2\cdots P_T
  =\cdots
  =\pi_\beta.
\end{align*}
}{{\revision The Gibbs transition probability for a sample
    corresponding to the update of each variable $\mbf u_t[k]$,
    $t\in\scr{T}$ satisfies the detailed balance
    equation~\eqref{eq:detailed-balance} for the desired distribution
    $p(u;\beta)$~\cite{Brooks1998}, which means that we also have global
    balance.} Detailed balance is not preserved through the full
  variable sweep, but global balance
  is\ifthenelse{\boolean{techrep}}{, which means that the state
    transition matrix associated with the samples
    $\mbf u[1],\mbf u[K+1], \dots$ has the desired invariant
    distribution.}{.}}
\ifthenelse{\boolean{techrep}}{

\subsubsection{Regularity and Convergence}

}{}The condition
\begin{align}\label{eq:positive-desired}
  p(u;\beta)>0, \quad\forall u\in\scr{U},
\end{align}
guarantees that for every two possible values of the Markov chain's
state, there is one path of nonzero probability that takes the state
from one value to the other over a full variable sweep (essentially by
changing one variable at a time). This means that the Markov chain
generated by Gibbs sampling is regular over the $T$ steps of a
variable sweep. \ifthenelse{\boolean{techrep}}{ To verify that this is
  so, imagine that at the start of a sweep ($t=1$) we have some
  arbitrary $\mbf u[k]=u\in\scr{U}$ and we want to verify that there is
  nonzero probability that we will end up at some other arbitrary
  $\mbf u[k+n]=u'\in\scr{U}$ at the end of the sweep.

\smallskip

In the first sweep step ($t=1$) only the first variable of $\mbf u[k]$
changes, but \eqref{eq:positive-desired} guarantees that there is a
positive probability that $\mbf u[k]$ we will transition to some
$\mbf u[k+1]$ for which the first variable matches that of $u'$, i.e.,
$\mbf u_1[k+1]=u'_1$, while all other variables remain unchanged. At
the second sweep step ($t=2$) only the first variable of $\mbf u[k+1]$
changes and now \eqref{eq:positive-desired} guarantees that there is a
positive probability that the second variable will also match that of
$u'$, i.e., we will get $\mbf u_1[k+1]=u'_1$ and
$\mbf u_2[k+1]=u'_2$. Continuing this reasoning, we conclude that
there is a positive probability that by the last sweep step the whole
$\mbf u[k+T]$ matches $u'$.

\smallskip

Having concluded that Gibbs sampling generates a regular Markov chain
across a full variable sweep,
Theorem~\ref{th:fundamental-Markov-chains} allow us to conclude that
the distribution of $\mbf u[k]$ converges exponentially fast to the
desired distribution. However, we shall see shortly that the mixing
time can scale poorly with the temperature parameter $\beta$.
}{}

\ifthenelse{\boolean{techrep}}{
\subsubsection{The binary case}

For a cost function $J(u)$ that takes binary values in $\{0,1\}$ and
each $u\in\scr{U}$ is an $n$-tuple of binary values one can compute
bounds on the values of the transition probabilities
$p_\mrm{sweep}(u' | u)$, $u,u' \in\scr{U}$ over the $n$ steps of a
variable sweep:
\begin{align*}
  p_\mrm{sweep}(u' | u)\ge
  \begin{cases}
    \frac{1}{1+e^{-\beta}}\Big(\frac{e^{-\beta}}{1+e^{-\beta}}\Big)^{n-1}&J(u')=0\\
    \Big(\frac{e^{-\beta}}{1+e^{-\beta}}\Big)^n&J(u')=1.
  \end{cases}
\end{align*}
from which we can conclude that
\begin{align*}
  p_{\min}
  &\eqdef \sum_{u' \in\scr{U}} \min_{u\in\scr{U}} p_\mrm{sweep}(u'|u)\\
  &\ge \Big(\frac{e^{-\beta}}{1+e^{-\beta}}\Big)^{n-1}
  \frac{|\scr{U}^*|+e^{-\beta}(2^n-|\scr{U}^*|)}{1+e^{-\beta}}\\
  p_{\max}&\eqdef \sum_{u' \in\scr{U}} \max_{u\in\scr{U}} p_\mrm{sweep}(u'|u)
  \le 2^n,
\end{align*}
where $\scr{U}^*\eqdef\{u\in\scr{U}:J(u)=0\}$. Using this in the
bound~\eqref{eq:Hoffman}, we get
\begin{align*}
  \frac{1}{\log\lambda^{-1}}
  \le \frac{
    (e^\beta+1)^n
    +(e^\beta-1)\frac{|\scr{U}^*|}{2^n}
    +1
    }
  {4(1-e^{-\beta})|\scr{U}^*|},
\end{align*}
which shows a bound for the mixing time growing exponentially with
$T\beta$. This bound is typically very conservative and, in fact, when
$\scr{U}^*$ has a single element, it is possible to obtain the
following much tighter bounds:
\begin{align*}
  p_{\min}
  &\eqdef \sum_{u'\in\scr{U}} \min_{u\in\scr{U}} p_\mrm{sweep}(u' | u)
  \ge\frac{e^{-\beta}(2^n-1)+1}{2^{n-1}(1+e^{-\beta})}\\
  p_{\max}
  &\eqdef \sum_{u'\in\scr{U}} \max_{u\in\scr{U}} p_\mrm{sweep}(u' | u)
  \le 1+\frac{1}{1+e^{-\beta}} \le 2,
\end{align*}
that lead to
\begin{align*}
  \frac{1}{\log\lambda^{-1}}
  \le\frac{1+2e^{-\beta}+2^{-n}(1-e^{-\beta})}{e^{-\beta}+2^{-n}(1-e^{-\beta})}
  \xrightarrow[n\to\infty]{}e^\beta+2,
\end{align*}
showing that the mixing time remains bounded as $T \to\infty$, but increases
exponentially as $\beta$ increases (temperature decreases).
}{}

\subsection{Parallel tempering}

Tempering 
decreases the mixing time of a Markov chain by
creating high-probability ``shortcuts'' between states.
It is applicable whenever we can embed the desired distribution into a
family of distributions parameterized by a parameter
$\beta\in[\beta_{\min},\beta_{\max}]$, with the property that we have slow mixing
for the desired distribution, which corresponds to
$\beta=\beta_{\max}$, but we have fast mixing for the distribution
corresponding to $\beta=\beta_{\min}$; which is typically for the Boltzmann
distribution~\eqref{eq:Boltzmann}.
The key idea behind tempering is then to select $M$ values
\begin{align*}
  \scr{B}\eqdef \{\beta_1\eqdef \beta_{\min} < \beta_2 <\beta_3 < \cdots < \beta_M\eqdef \beta_{\max} \}
\end{align*}
and generate samples from the joint distribution
\begin{align}\label{eq:tempering-desired-joint}
  p(u^\beta:\beta\in\scr{B})\eqdef \prod_{\beta\in\scr{B}} p(u^\beta;\beta)
\end{align}
that corresponds to $M$ independent random variablea $\mbf u^\beta$, one
for each parameter value $\beta\in\scr{B}$. We group these variables as an
$M$-tuple and denote the joint Markov chain by
\begin{align*}
  \mbf u[k] \eqdef (\mbf u^\beta[k]: \beta\in\scr{B}),
\end{align*}
Eventually, from each $M$-tuple we only use the samples
$\mbf u^\beta[k]$, $\beta=\beta_{\max}$ that correspond to the desired
distribution.

\subsubsection{General algorithm}

Tempering can be applied to any MCMC method associated with a regular
Markov chain with transition probabilities $p(u'| u;\beta)$ and strictly
positive transition matrices $P_\beta$, $\beta\in\scr{B}$.
The algorithm uses a \emph{flip function}
defined by
\begin{align}\label{eq:tempering-flip}
  f_\mrm{flip}(u^{\beta_j},u^{\beta_{j+1}})
  &=\min\Big\{
  \frac{p(u^{\beta_{j+1}};\beta_j)}{p(u^{\beta_{j+1}};\beta_{j+1})}
  \frac{p(u^{\beta_j};\beta_{j+1})}{p(u^{\beta_j};\beta_j)}
  ,1\Big\},
\end{align}
and operates as follows:

\smallskip

\begin{algorithm}[Tempering]~\small
  \label{alg:tempering}
\begin{senumerate}{0.5ex}
\item Pick arbitrary $\mbf u[1]=(\mbf u^\beta[1]\in\scr{U}: \beta\in\scr{B})$.
\item\label{en:loop-tempering} Set $k=1$ and repeat until enough samples are collected:
  \begin{enumerate}
  \item\label{en:tempering-parallel} For each $\beta\in\scr{B}$, sample
    $\mbf u^\beta[k+1]$ with probability
    $p(u'| \mbf u^\beta[k];\beta)$ and increment $k$.

  \item\label{en:tempering-all-flips} \emph{Tempering sweep}: For each $j\in\{1,\dots,M-1\}$:

    \begin{itemize}
    \item\label{en:tempering-flip}  Compute the flip probability
      \begin{align}\label{eq:p-flip}
        p_\mrm{flip}=f_\mrm{flip}(\mbf u^{\beta_j}[k],\mbf u^{\beta_{j+1}}[k])
      \end{align}
      and set
      \begin{align*}
        \mbf u[k+1] =
        \begin{cases}
          \mbf {\tilde u}[k]
          & \text{with prob.~} p_\mrm{flip},\\
          \mbf u[k] & \text{with prob.~} 1-p_\mrm{flip},
        \end{cases}
      \end{align*}
      where $\mbf {\tilde u}[k]$ is a version of $\mbf u[k]$ with the
      entries $\mbf u^{\beta_j}[k]$ and $\mbf u^{\beta_{j+1}}[k]$ flipped; and
      increment $k$.\frqed
    \end{itemize}
  \end{enumerate}
\end{senumerate}
\end{algorithm}

\smallskip

Step~\ref{en:tempering-parallel} corresponds to one step
of a base MCMC algorithm (e.g., Gibbs sampling), for each value of
$\beta\in\scr{B}$.  For the Boltzmann distribution~\eqref{eq:Boltzmann}, the
flip function is given by
\begin{align*}
  f_\mrm{flip}(u^{\beta_j},u^{\beta_{j+1}})
  &=\min\Big\{e^{-(\beta_j-\beta_{j+1})(J(u^{\beta_{j+1}})-J(u^{\beta_j}))},1\Big\},
\end{align*}
which means that the variables $\mbf u^{\beta_j}[k]$ and
$\mbf u^{\beta_{j+1}}[k]$ are flipped with probability one whenever
$J(\mbf u^{\beta_j})< J(\mbf u^{\beta_{j+1}})$. Intuitively, the tempering
sweep in step \ref{en:tempering-all-flips} quickly brings to
$\mbf u^{\beta_{\max}}[k]$ low-energy/low-cost samples that may have been
``discovered'' by other $\mbf u^\beta[k]$, $\beta<\beta_{\max}$ with better
mixing.

\subsubsection{Balance}

{\revision Since the sample extractions in step~\ref{en:tempering-parallel} are
independent, the transition probability corresponding to this
step is given by
\begin{align*}
  p\big( ({u'}^\beta:\beta\in\scr{B}) \,|\, (u^\beta:\beta\in\scr{B}) \big)
  = \prod_{\beta\in\scr{B}} p({u'}^\beta|u^\beta;\beta),
\end{align*}
which satisfies the global balance equation for the joint distribution
in~\eqref{eq:tempering-desired-joint}.
\ifthenelse{\boolean{techrep}}{Picking some $j\in\{1,\dots,M-1\}$, the
  transition probability for step~\ref{en:tempering-flip} is given by
\begin{multline*}
  p\Big( ({u'}^\beta:\beta\in\scr{B}) \,|\, (u^\beta:\beta\in\scr{B}) \Big)=\\
  \begin{cases}
    \scriptstyle
    f_\mrm{flip}(u^{\beta_j},u^{\beta_{j+1}})
    &
    \scriptstyle
    {u'}^{\beta_j}=u^{\beta_{j+1}},\, u^{\beta_j}={u'}^{\beta_{j+1}},\,
    {u'}^\beta=u^\beta,\, 
    \beta\not\in \{ \beta_j,\beta_{j+1}\}\\
    0 & \text{otherwise}
    \ifthenelse{\boolean{techrep}}{,}{.}
  \end{cases}
\end{multline*}
and therefore
\begin{multline*}
  p\Big( ({u'}^\beta:\beta\in\scr{B}) \,|\, (u^\beta:\beta\in\scr{B}) \Big)
  p(u^\beta:\beta\in\scr{B})\\
  =\begin{cases}
    \scriptstyle
    p(u^{\beta_j};\beta_j)p(u^{\beta_{j+1}};\beta_{j+1})
    p(u^\beta:\beta\not\in \{ \beta_j,\beta_{j+1}\})
    f_\mrm{flip}(u^{\beta_j},u^{\beta_{j+1}})
    \\\quad\text{if~}
    \scriptstyle
    {u'}^{\beta_j}=u^{\beta_{j+1}},\; u^{\beta_j}={u'}^{\beta_{j+1}},
    {u'}^\beta=u^\beta, \;\forall \beta\not\in \{ \beta_j,\beta_{j+1}\},\\
    0
    \quad
    \text{otherwise},
  \end{cases}
\end{multline*}
whereas
\begin{multline*}
  p\Big( (u^\beta:\beta\in\scr{B}) \,|\, ({u'}^\beta:\beta\in\scr{B}) \Big)
  p({u'}^\beta:\beta\in\scr{B})\\
  =\begin{cases}
    \scriptstyle
    p({u'}^{\beta_j};\beta_j)p({u'}^{\beta_{j+1}};\beta_{j+1})
    p({u'}^\beta:\beta\not\in \{ \beta_j,\beta_{j+1}\})
    f_\mrm{flip}({u'}^{\beta_j},{u'}^{\beta_{j+1}})
    \\\quad\text{if~}
    \scriptstyle
    u^{\beta_j}={u'}^{\beta_{j+1}},\; {u'}^{\beta_j}=u^{\beta_{j+1}},
    {u'}^\beta=u^\beta, \;\forall \beta\not\in \{ \beta_j,\beta_{j+1}\},\\
    0
    \quad
    \text{otherwise}
  \end{cases}\\
  =\begin{cases}
    \scriptstyle
    p(u^{\beta_{j+1}};\beta_j)p(u^{\beta_j};\beta_{j+1})
    p(u^\beta:\beta\not\in \{ \beta_j,\beta_{j+1}\})
    f_\mrm{flip}(u^{\beta_{j+1}},u^{\beta_j})
    \\\quad\text{if~}
    \scriptstyle
    u^{\beta_j}={u'}^{\beta_{j+1}},\; {u'}^{\beta_j}=u^{\beta_{j+1}},
    {u'}^\beta=u^\beta, \;\forall \beta\not\in \{ \beta_j,\beta_{j+1}\},\\
    0
    \quad
    \text{otherwise}.
  \end{cases}
\end{multline*}
To get detailed balance, we thus need
\begin{multline*}
  p(u^{\beta_j};\beta_j)p(u^{\beta_{j+1}};\beta_{j+1})
  f_\mrm{flip}(u^{\beta_j},u^{\beta_{j+1}})\\
  =
  p(u^{\beta_{j+1}};\beta_j)p(u^{\beta_j};\beta_{j+1})
  f_\mrm{flip}(u^{\beta_{j+1}},u^{\beta_j}).
\end{multline*}
This equality always holds for the flip function defined
in~\eqref{eq:tempering-flip}, for which we always have either
$f_\mrm{flip}(u^{\beta_j},u^{\beta_{j+1}})=1$ or
$f_\mrm{flip}(u^{\beta_{j+1}},u^{\beta_j})=1$. This guarantees detailed
balance for all $M-1$ steps in~\ref{en:tempering-flip} and therefore
global balance for all the combined steps
in~\ref{en:tempering-parallel} and \ref{en:tempering-all-flips} for a
full tempering sweep.}{The flip function in~\eqref{eq:tempering-flip}
guarantees detailed balance for all $M-1$ steps
in~\ref{en:tempering-flip} and therefore global balance for all the
combined steps in~\ref{en:tempering-parallel} and
\ref{en:tempering-all-flips}, within a full tempering sweep~\cite{HespanhaMar2024b_arXiv}.
}}

\subsubsection{Regularity and Convergence}

Assuming that for each $\beta\in\scr{B}$, the Markov chain generated by
$p(u'| u;\beta)$ is regular with a strictly positive transition matrix
$P_\beta$, any possible combination of states
$\mbf u[k] = (u^\beta[k]: \beta\in\scr{B})$ at time $k$ can transition in one
time step to any possible combination of states
$\mbf u[k+1] = (u^\beta[k]: \beta\in\scr{B})$ at time $k+1$ at
step~\ref{en:tempering-parallel}. This means that this step
corresponds to a transition matrix $P$ with strictly positive
entries. In contrast, each flip step corresponds to a transition
matrix $P_{\mrm{flip } j}$ that is non-negative and right-stochastic,
but typically has many zero entries. However, each matrix $P_{\mrm{flip } j}$
cannot have any row that is identically zero (because rows must add up
to 1). This suffices to conclude that any product of the form
\begin{align*}
  Q=P_{\mrm{flip } 1}P_{\mrm{flip } 2}\cdots P_{\mrm{flip } M-1} P
\end{align*}
must have all entries strictly positive. This transition matrix
corresponds to the transition from the start of step
\ref{en:tempering-all-flips} in one ``tempering sweep'' to the start
of the same step at the next sweep and defines a regular Markov
chain. {\revision Theorem~\ref{th:fundamental-Markov-chains} thus
  allow us to conclude that the distribution at the start of
  step~\ref{en:tempering-all-flips} converges to the desired invariant
  distribution. Since every step satisfies global balance, the
  invariant distribution is preserved across every step, so the
  distributions after every step also converges to the invariant
  distribution.}


\section{MCMC FOR OPTIMAL CONTROL}\label{se:MCMC-Koopman}

Our goal is to optimize a
criterion of the form~\eqref{eq:Koopman-criterion}.
In the context of MCMC sampling from the Boltzmann distribution
\eqref{eq:Boltzmann}, this corresponds to a Gibbs update for the
control $\mbf u^\beta_t[k+1]$ with a distribution~\eqref{eq:Gibbs-i} that
can be computed using
\begin{align*}
  \frac{e^{-\beta J(u_t,\mbf u^\beta_{-t}[k])}}{\sum_{\bar u_t} e^{-\beta J(\bar u_t,\mbf u^\beta_{-t}[k])}}
  =\frac{e^{-\beta \mbf c^\beta_t[k] A(u_t)\mbf x^\beta_t[k]}}{
    \sum_{\bar u_t} e^{-\beta e^{-\beta \mbf c^\beta_t[k] A(\bar u_t)\mbf x^\beta_t[k]}}},
\end{align*}
and a tempering flip probability in \eqref{eq:p-flip} computed using
\begin{align*}
  p_\mrm{flip}
  \ifthenelse{\boolean{techrep}}{
    &=\min\Big\{e^{-(\beta_j-\beta_{j+1})(J(\mbf u^{\beta_{j+1}}[k])-J(\mbf u^{\beta_j}[k]))},1\Big\}\\
    }{}
  &=\min\Big\{e^{-(\beta_j-\beta_{j+1})c(\mbf x^{\beta_{j+1}}_T[k])-\mbf x^{\beta_j}_T[k]))},1\Big\},
\end{align*}
where
\begin{align}
  \label{eq:c-k}
  \mbf c^\beta_t[k]&\eqdef c A(\mbf u^\beta_T[k])\cdots A(\mbf u^\beta_{t+1}[k]), \\
  \label{eq:x-k}
  \mbf x^\beta_t[k]&\eqdef A(\mbf u^\beta_{t-1}[k])\cdots A(\mbf u^\beta_1[k])x_1.
\end{align}
The following algorithm implements Gibbs sampling with tempering using
recursive formulas to evaluate~\eqref{eq:c-k}--\eqref{eq:x-k}.

\smallskip

\begin{algorithm}[Tempering for Koopman optimal control]~\small
  \label{alg:tempering-Koopman}
\begin{senumerate}{0.5ex}
\item Pick arbitrary
  $\mbf u[1]=(\mbf u_t^\beta[1]\in\scr{U}: t\in\scr{T}, \beta\in\scr{B})$.
\item\label{en:loop-Gibbs-tempering} Set $k=1$ and repeat until enough
  samples are collected:
  \begin{enumerate}
  \item \emph{Gibbs sampling}: For each $\beta\in\scr{B}$:
    \begin{itemize}
    \item \revision Set $\mbf x^\beta_1[k]=x_1$, $\mbf c^\beta_T[k]=c$ and, $\forall t\in\scr{T}$,
      \begin{align} \label{eq:update-C}
        &\mbf c^\beta_t[k]= \mbf c^\beta_{t+1}[k] A(\mbf u^\beta_{t+1}[k]).
      \end{align}
    \item \emph{Variable sweep}: For each $t\in\scr{T}$
      \begin{itemize}
      \item Sample $\mbf u^\beta_t[k+1]$ with distribution
        \begin{align}\label{eq:Gibbs-linear}
          \frac{e^{-\beta \mbf c^\beta_t[k] A(u_t)\mbf x^\beta_t[k]}}{
            \sum_{\bar u_t} e^{-\beta e^{-\beta \mbf c^\beta_t[k] A(\bar u_t)\mbf x^\beta_t[k]}}},
        \end{align}
      \item Set $\mbf u^\beta_{-t}[k+1]=\mbf u^\beta_{-t}[k]$, update
        \begin{align}
          \label{eq:update-X}
          &\mbf x^\beta_{t+1}[k+1]=A(\mbf u^\beta_t[k+1])\mbf x^\beta_t[k],
        \end{align}
        and increment $k$.
      \end{itemize}
    \end{itemize}

  \item \emph{Tempering sweep}: For each $j\in\{1,\dots,M-1\}$
    \begin{itemize}
    \item Compute the flip probability
      \begin{align*}
        p_\mrm{flip}=
        \min\Big\{e^{-(\beta_j-\beta_{j+1})c(\mbf x^{\beta_{j+1}}_T[k])-\mbf x^{\beta_j}_T[k]))},1\Big\},
      \end{align*}
      and set
      \begin{align*}
        \mbf u[k+1] =
        \begin{cases}
          \mbf {\tilde u}[k] &\text{with prob.~} p_\mrm{flip}\\
          \mbf u[k] & \text{with prob.~} 1-p_\mrm{flip},
        \end{cases}
      \end{align*}
      where $\mbf {\tilde u}[k]$ is a version of $\mbf u[k]$ with the
      entries $\mbf u^{\beta_j}[k]$ and $\mbf u^{\beta_{j+1}}[k]$ flipped; and
      increment $k$. \frqed
    \end{itemize}
  \end{enumerate}
\end{senumerate}
\end{algorithm}

\smallskip

\subsubsection*{Computation complexity and parallelization}

The bulk of the computation required by
Algorithm~\ref{alg:tempering-Koopman} lies in the computation of the
matrix-vector products that appear in \eqref{eq:update-C},
\eqref{eq:Gibbs-linear}, and \eqref{eq:update-X}, each of these
products requiring $O(n_\psi^2)$ floating-point operation for a total of
$O(MT(2+|\scr{U}|)n_\psi^2)$ operations.  The tempering step, does not
use any additional vector-matrix multiplications.
In contrast, a naif implementation of Gibbs sampling with tempering
would require $MT|\scr{U}|$ evaluations of the cost
function~\eqref{eq:Koopman-criterion}, each with computational complexity
$O(Tn_\psi^2)$. {\revision The reduction in computation complexity by a
  factor of $T|\scr{U}|/(2+|\scr{U}|)$ is especially significant when
  the time horizon is large. The price paid for the computational
  savings is that we need to store the
  vectors~\eqref{eq:update-C}, \eqref{eq:update-X}, with memory complexity
  $O(MTn_\psi)$.}

\smallskip

The computations of the matrix-vector products mentioned above are
independent across different values of $\beta\in\scr{B}$ and can be
performed in parallel. This means that tempering across $M$
temperatures can be computationally very cheap {\revision if
  enough computational cores and memory are available.}
\ifthenelse{\boolean{techrep}}{ In contrast, parallelization across a
  Gibbs variable sweep is not so easily parallelizable because, within
  one variable sweep, the value of each sample $\mbf u^\beta_t[k+1]$,
  typically depends on the values of previous samples
  $\mbf u^\beta_t[k]$ (for the same value of $\beta$). Tempering thus both
  decreases the mixing time of the Markov chain and opens the door for
  a high degree of parallelization.}{}

\section{Numerical Example}\label{se:numerical}


We tested the algorithm proposed in this paper on a Koopman model for
the Atari 2600 \emph{Assault} video game. The goal of the game is to
protect earth from small attack vessels deployed by an alien
mothership. The mother ship and attack {\revision vessels shoot} at the player's
ship and the player uses a joystick to dodge the incoming fire and
fire back at the alien ships. The player's ship is destroyed either if
it is hit by enemy fire or if it ``overheats'' by shooting at the
aliens. The player earns points by destroying enemy ships.

\smallskip

We used a Koopman model with $|\scr{U}|=4$ control action, which
correspond to ``move left'', ``move right'', ``shoot up'', and do
nothing. The optimization minimizes the  cost
\begin{align*}
  J(u)\eqdef
  \begin{cases}
    -\frac{1}{T} \sum_{t=1}^T r_t & \text{if player's ship not destroyed}\\
    +10-\frac{1}{T} \sum_{t=1}^T r_t & \text{if player's ship destroyed}
  \end{cases}
\end{align*}
where $r_t$ denotes the points earned by the player at time $t$.
This cost balances the tradeoff between taking some risk to collect
points to decrease $-\frac{1}{T} \sum_{t=1}^T r_t$, while not getting
destroyed and incurring the $+10$ penalty. {\revision For game play,
  this optimization is solved at every time step with a receding
  horizon\ifthenelse{\boolean{techrep}}{ starting at the current time
    $t$ and ending at time $t+T$, with only the first control action
    executed}{}. However, because the focus of this paper is on the
  solution
  to~\eqref{eq:nonlinear-process}--\eqref{eq:nonlinear-minimization},
  we present results for a single optimization starting from a typical
  initial condition. A time horizon $T=40$ was used in this section,
  which corresponds to a total number of control options
  $|\scr{U}|^T\approx1.2\times 10^{24}$.}

\smallskip

The state of the system is built directly from screen pixel
information. Specifically, the pixels are segmented into 5 categories
corresponding to the player's own ship, the player's horizontal fire,
the player's vertical fire, the attacker's ships and fire, and the
temperature bar. For the player's own ship and the attacker's
ships/fire, we consider pixel information from the current and last
screenshot, so that we have ``velocity'' information. {\revision The
  pixels of each category are used to construct ``spatial densities''
  using the entity-based approach described
  in~\cite{BlischkeHespanhaDec2023}, with observables of the form
\begin{align}\label{eq:smooth-densities}
  \varphi_{\ell j}(x)\eqdef \sum_i e^{-\lambda \|p_{\ell i}-c_{\ell j}\|^2},
\end{align}
where the index $\ell$ ranges over the 5 categories above, the summation
is taken over the pixels $p_{\ell i}$ associated with category $\ell$, and
the $c_{\ell j}$ are fixed 
points in the screen. The densities associated with the 5 categories
are represented by 50, 50, 100, 200, 16 points $c_{\ell j}$,
respectively. The observables in \eqref{eq:smooth-densities}, together
with the value of the optimization criterion,} are used to form a
lifted state $\psi_t$ 
with dimension $n_\psi=667$\ifthenelse{\boolean{techrep}}{ (see
  Table~\ref{tab:assault-state})}{}. {\revision The matrices $A(u)$
  in~\eqref{eq:Koopman-criterion} were estimated using 500 game traces
  using random inputs. Collecting data from the Atari simulator and
  lifting the state took about 1h45min, whereas solving the
  least-squares problems needed to obtain the Koopman matrices took
  less than 1 sec.}

\smallskip

\ifthenelse{\boolean{techrep}}{
  \begin{table}
    \caption{State of the Koopman model}
    \label{tab:assault-state}
    \centering
    \begin{tabular}{l|c}
      Description&dimension\\\hline
      player's ship current pixels & 50 \\
      player's ship last pixels & 50 \\
      player's horizontal fire current pixels & 50\\
      player's vertical fire current pixels & 100\\
      attacker's ships current pixels & 200 \\
      attacker's ships last pixels & 200 \\
      temperature bar current pixels & 16\\
      optimization criterion & 1\\\hline
      Koopman state dimension $n_\psi$& 667
    \end{tabular}
\end{table}
}

\begin{figure}[tbh]
  \parbox{\columnwidth}{
    \centering
    \includegraphics[width=.98\columnwidth]{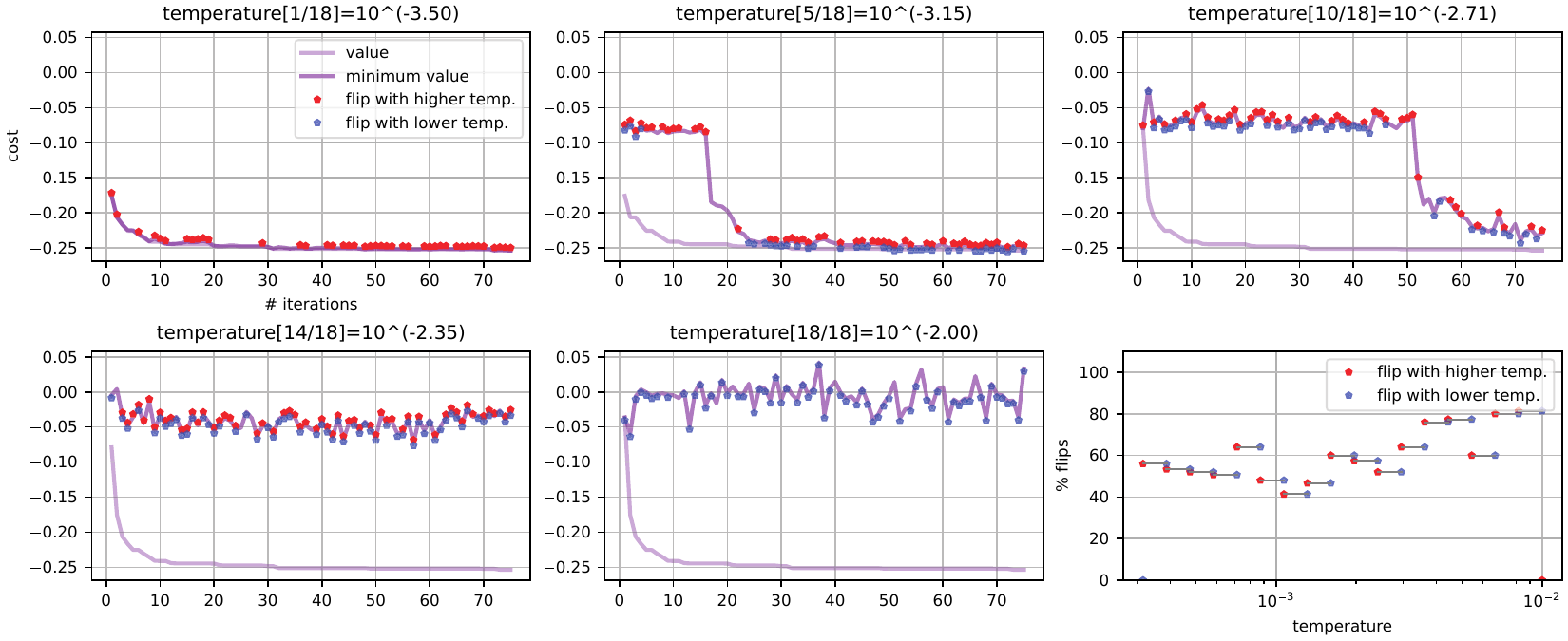}
  }
  \caption{\revision Typical run of
    Algorithm~\ref{alg:tempering-Koopman}: The first 5 plots show the
    cost ($y$-axis) at the end of each iteration ($x$-axis), as well
    as the minimum cost found so far, for a specific temperature
    $\beta$. For this run, 12 logarithmic scaled temperatures were used,
    but only 5 are shown here. In each plot, red and blue dots mark
    iterations where flips occurred during the tempering sweep. The
    bottom-right plot shows the total number of flips across
    ``adjacent'' temperatures (in the $x$-axis), as a percentage of
    the total number of iterations (in the
    $y$-axis).\label{fig:typical-run}}
\end{figure}

Figure~\ref{fig:typical-run} depicts a typical run of
Algorithm~\ref{alg:tempering-Koopman}, showing flipping of samples
across adjacent temperatures in 40-80\% of the tempering sweeps.
{\revision In the remainder of this section we compare the performance
  of Algorithm~\ref{alg:tempering-Koopman} with several alternatives.
  All run times refer to Julia implementations on a 2018
  MacBook Pro with a 2.6GHz 
  Intel Core
  i7 CPU.
}

\smallskip

\begin{figure}[tbh]
  \centering
    \parbox{\columnwidth}{
      \includegraphics[width=.98\columnwidth]{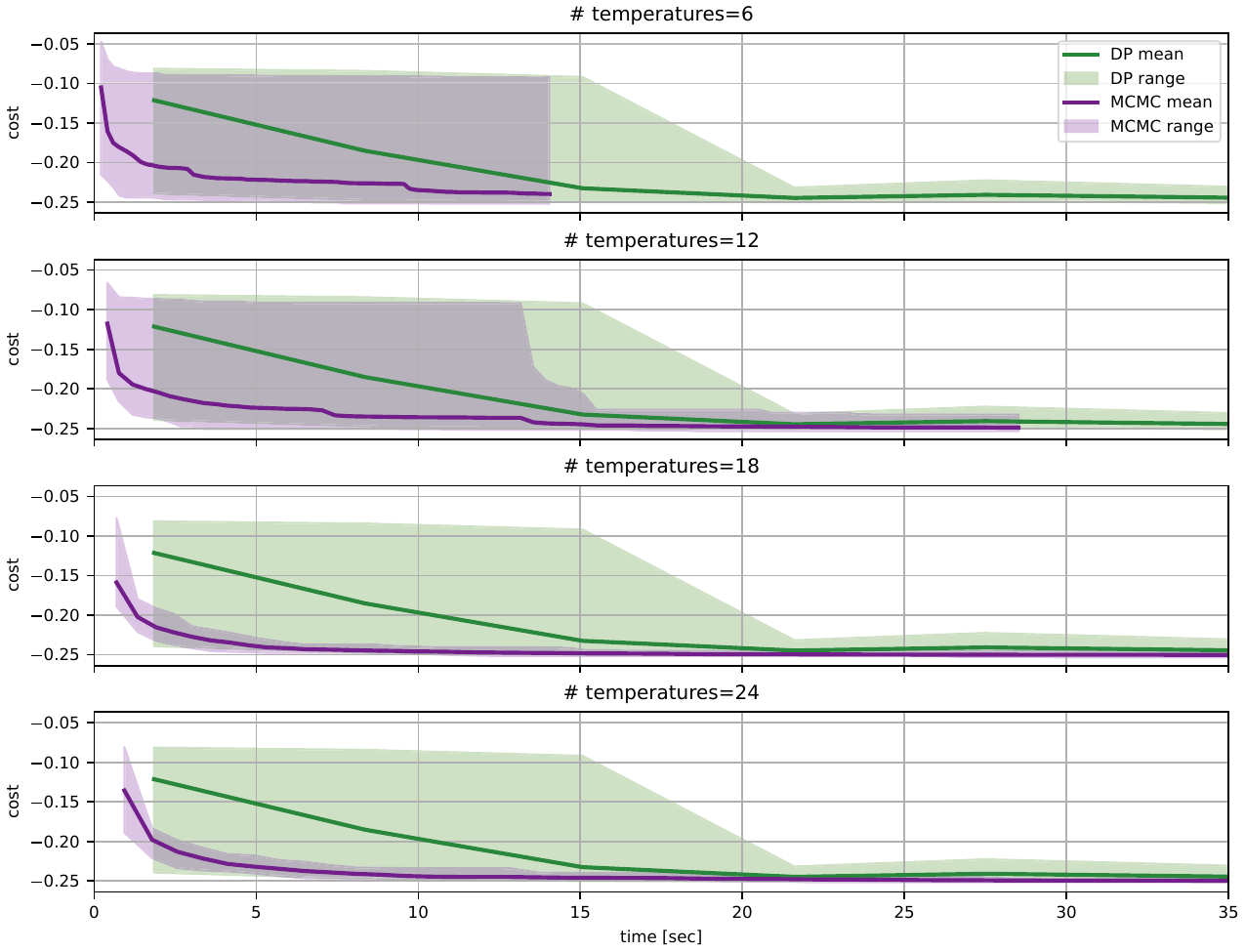}
    }
    \caption{Run-time comparison between
      Algorithm~\ref{alg:tempering-Koopman} (purple) and the
      dynamic programming algorithm in~\cite{BlischkeHespanhaDec2023}
      (green). The $x$-axis denotes run-time and the $y$-axis the cost
      achieved, with the solid lines showing the average cost across
      15 different runs and the shaded areas the whole range of costs
      obtained over those runs. The different plots correspond to
      different values for the number of
      temperatures $M$. \label{fig:DP-BP-comparison}}
\end{figure}

Figure~\ref{fig:DP-BP-comparison} compares
Algorithm~\ref{alg:tempering-Koopman} with the algorithm
in~\cite{BlischkeHespanhaDec2023}, which exploits the piecewise-linear
structure of the cost-to-go to efficiently represent and evaluate the
value function and also to dynamically prune the search tree.  Due to
the need for exploration, this algorithm ``protects'' from pruning a
random fraction of tree-branches (see~\cite{BlischkeHespanhaDec2023}
for details).  Both algorithms used 6 CPU cores. For
Algorithm~\ref{alg:tempering-Koopman}, each core executed one sweep of
the tempering algorithm and for the algorithm
in~\cite{BlischkeHespanhaDec2023}, the 6 cores were used by BLAS to
speedup matrix multiplication.
%
%
Both algorithms were executed multiple times and the plots show the
costs obtained as a function of run time. For
Algorithm~\ref{alg:tempering-Koopman}, the run time is
{\revision controlled} by the number of samples drawn. For the
algorithm in~\cite{BlischkeHespanhaDec2023}, the run time is
{\revision controlled} by the number of vectors used to represent the
value function. Both algorithms eventually discover comparable
``optimal'' solutions, but Algorithm~\ref{alg:tempering-Koopman}
consistently finds a lower cost faster.



\begin{figure}[htp]
  \centering
    \parbox{\columnwidth}{
      \includegraphics[width=.98\columnwidth]{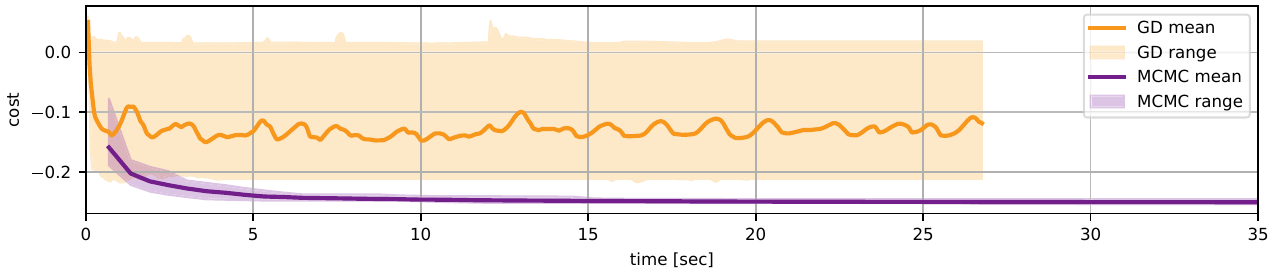}
    }
    \\
    \parbox{\columnwidth}{
    \includegraphics[width=.98\columnwidth]{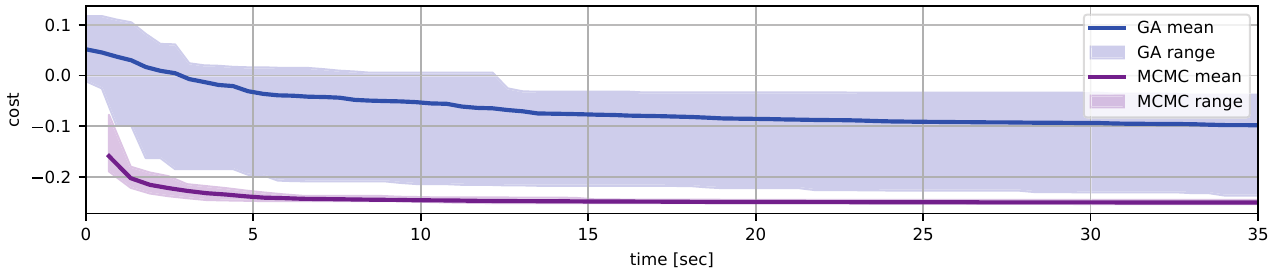}
    }
    \caption{\revision Run-time comparison between
      Algorithm~\ref{alg:tempering-Koopman} with 18 temperatures
      (purple) and a gradient descent algorithm (top, orange) and a
      genetic algorithm (bottom, blue).  The meaning of the solid
      lines and shaded areas is the same as in
      Figure~\ref{fig:DP-BP-comparison}. \label{fig:GD-GA-BP-comparison}}
\end{figure}

\smallskip

{\revision The \revision top plot in Figure~\ref{fig:GD-GA-BP-comparison} compares
Algorithm~\ref{alg:tempering-Koopman} with a gradient descent solver
that minimizes the following continuous relaxation of the
cost~\eqref{eq:Koopman-criterion}:
\begin{align}\label{eq:optimization-relaxed}
  &J_\mrm{relax}(\mu)\eqdef c\psi_T, &
  & \psi_{t+1} = \Big(\sum_{u\in\scr{U}} \mu_u(t) A(u)\Big) \psi_t
\end{align}
where the optimization variables $\mu_u(t)\in[0,1]$ are required to
satisfy $\sum_{u\in\scr{U}} \mu_u(t)=1$, $\forall t\in\scr{T}$.  The global minimum
$J_\mrm{relax}^*$ of~\eqref{eq:optimization-relaxed} would match that
of~\eqref{eq:Koopman-criterion}, if we forced the $\mu_u(t)$ to take
binary values in $\{0,1\}$. Otherwise, $J_\mrm{relax}^*$ provides a
lower bound for~\eqref{eq:optimization-discrete}. We
minimized~\eqref{eq:optimization-relaxed} using the
toolbox~\cite{GradDescent.jl}. The results shown were obtained using
Nesterov-accelerated adaptive moment estimation~\cite{GradDescent.jl},
with $\eta=0.5$, which resulted in the best performance among the
algorithms supported by~\cite{GradDescent.jl}. The constraints on
$\mu_u(t)$ were enforced by minimum-distance projection into the
constraint set. In general, gradient descent converges quickly, but to
a local minima of~\eqref{eq:optimization-relaxed} with cost higher
then the minimum found by Algorithm~\ref{alg:tempering-Koopman}
for~\eqref{eq:optimization-discrete}; in spite of the fact that the
global minima of~\eqref{eq:optimization-relaxed} is potentially
smaller than that of~\eqref{eq:optimization-discrete}.

\smallskip

The bottom plot in Figure~\ref{fig:GD-GA-BP-comparison} compares
Algorithm~\ref{alg:tempering-Koopman} with a genetic optimization
algorithm that also resorts to stochastic exploration. This type of
algorithm simulates a ``population'' of candidate solutions to the
optimization~\eqref{eq:optimization-discrete}, which evolves by
mutation, crossover, and selection.
We minimized~\eqref{eq:optimization-relaxed} using the
toolbox~\cite{Evolutionary.jl}. The results in
Figure~\ref{fig:GD-GA-BP-comparison} used a population size $P=50$,
selection based on uniform ranking (which selects the best $\mu=2$
individuals with probability $1/\mu$), binary crossover (which combines
the solution of the two parent with a single crossover point), a
mutation rate of 10\%, and a probability of mutation for each ``gene''
of 5\%. These parameters resulted in the lowest costs we could achieve
among the options provided by~\cite{Evolutionary.jl}, but still
significantly higher than the costs obtained with
Algorithm~\ref{alg:tempering-Koopman}.
}







\bibliographystyle{abbrvnat}
\bibliography{Kmcmc-TR-bib}

\begin{thebibliography}{31}
\providecommand{\natexlab}[1]{#1}
\providecommand{\url}[1]{\texttt{#1}}
\expandafter\ifx\csname urlstyle\endcsname\relax
  \providecommand{\doi}[1]{doi: #1}\else
  \providecommand{\doi}{doi: \begingroup \urlstyle{rm}\Url}\fi

\bibitem[Evo()]{Evolutionary.jl}
{\footnotesize \url{https://github.com/wildart/Evolutionary.jl}}.

\bibitem[Gra()]{GradDescent.jl}
{\footnotesize \url{https://github.com/jacobcvt12/GradDescent.jl}}.

\bibitem[Aadit et~al.(2022)Aadit, Grimaldi, Carpentieri, Theogarajan, Martinis,
  Finocchio, and Camsari]{aadit2022massively}
N.~Aadit, A.~Grimaldi, M.~Carpentieri, L.~Theogarajan, J.~Martinis,
  G.~Finocchio, and K.~Y. Camsari.
\newblock Massively parallel probabilistic computing with sparse {I}sing
  machines.
\newblock \emph{Nature Electronics}, 5\penalty0 (7):\penalty0 460--468, 2022.

\bibitem[Bakker et~al.(2019)Bakker, Nowak, and Rosenthal]{bakker2019learning}
C.~Bakker, K.~E. Nowak, and W.~S. Rosenthal.
\newblock Learning {K}oopman operators for systems with isolated critical
  points.
\newblock In \emph{Proc.~of the 58th IEEE Conf.~on Decision and Contr.}, pages
  7733--7739, 2019.

\bibitem[Bengea and DeCarlo(2005)]{bengea2005optimal}
S.~C. Bengea and R.~A. DeCarlo.
\newblock Optimal control of switching systems.
\newblock \emph{Automatica}, 41\penalty0 (1):\penalty0 11--27, 2005.

\bibitem[Bevanda et~al.(2021)Bevanda, Sosnowski, and
  Hirche]{BevandaSosnowskiHirche2021}
P.~Bevanda, S.~Sosnowski, and S.~Hirche.
\newblock {K}oopman operator dynamical models: Learning, analysis and control.
\newblock \emph{Annual Reviews in Control}, 52:\penalty0 197--212, 2021.

\bibitem[Blischke and Hespanha(2023)]{BlischkeHespanhaDec2023}
M.~Blischke and J.~P. Hespanha.
\newblock Learning switched {K}oopman models for control of entity-based
  systems.
\newblock In \emph{Proc.~of the 62th IEEE Conf.~on Decision and Contr.}, Dec.
  2023.

\bibitem[Brooks(1998)]{Brooks1998}
S.~Brooks.
\newblock Markov chain {M}onte {C}arlo method and its application.
\newblock \emph{J.~of the Royal Statistical Society: series D (the
  Statistician)}, 47\penalty0 (1):\penalty0 69--100, 1998.

\bibitem[Brunton et~al.(2016)Brunton, Brunton, Proctor, and
  Kutz]{BruntonBruntonProctorKutz2016}
S.~Brunton, B.~Brunton, J.~Proctor, and J.~N. Kutz.
\newblock {K}oopman invariant subspaces and finite linear representations of
  nonlinear dynamical systems for control.
\newblock \emph{PloS one}, 11\penalty0 (2), 2016.

\bibitem[{\v{C}}ern{\`y}(1985)]{Cerny1985thermodynamical}
V.~{\v{C}}ern{\`y}.
\newblock Thermodynamical approach to the traveling salesman problem: An
  efficient simulation algorithm.
\newblock \emph{J.~Opt.~Theory and Applications}, 45:\penalty0 41--51, 1985.

\bibitem[Earl and Deem(2005)]{Earl2005parallel}
D.~J. Earl and M.~W. Deem.
\newblock Parallel tempering: Theory, applications, and new perspectives.
\newblock \emph{Physical Chemistry Chemical Physics}, 7\penalty0 (23):\penalty0
  3910--3916, 2005.

\bibitem[Folkestad and Burdick(2021)]{FolkestadBurdick2021}
C.~Folkestad and J.~W. Burdick.
\newblock {K}oopman {NMPC}: {K}oopman-based learning and nonlinear model
  predictive control of control-affine systems.
\newblock In \emph{Proc.~of the IEEE Int.~Conf.~on Robot.~and Automat. (ICRA)},
  pages 7350--7356, 2021.

\bibitem[Gelfand and Smith(1990)]{Gelfand1990sampling}
A.~E. Gelfand and A.~F. Smith.
\newblock Sampling-based approaches to calculating marginal densities.
\newblock \emph{J.~of the American Statistical Assoc.}, 85\penalty0
  (410):\penalty0 398--409, 1990.

\bibitem[Geman and Geman(1984)]{GemanGeman1984stochastic}
S.~Geman and D.~Geman.
\newblock Stochastic relaxation, {G}ibbs distributions, and the {B}ayesian
  restoration of images.
\newblock \emph{IEEE Trans.~on Pattern Anal.~and Machine Intell.},
  PAMI-6\penalty0 (6):\penalty0 721--741, 1984.

\bibitem[Geyer(1991)]{Geyer1991markov}
C.~J. Geyer.
\newblock {M}arkov chain {M}onte {C}arlo maximum likelihood.
\newblock In \emph{Computing Science and Statistics: Proc.~of the 23rd Symp.~on
  the Interface}, pages 156--163, 1991.

\bibitem[Haseli and Cortés(2023)]{HaseliCortes2023}
M.~Haseli and J.~Cortés.
\newblock Modeling nonlinear control systems via {K}oopman control family:
  Universal forms and subspace invariance proximity, 2023.
\newblock arXiv: 2307.15368.

\bibitem[Hukushima and Nemoto(1996)]{Hukushima1996exchange}
K.~Hukushima and K.~Nemoto.
\newblock Exchange {M}onte {C}arlo method and application to spin glass
  simulations.
\newblock \emph{J.~of the Physical Society of Japan}, 65\penalty0 (6):\penalty0
  1604--1608, 1996.

\bibitem[Kemeny and Snell(1976)]{KemenySnell1976}
J.~G. Kemeny and J.~L. Snell.
\newblock \emph{Finite {M}arkov chains}.
\newblock Springer-Verlag, 1976.

\bibitem[Kirkpatrick et~al.(1983)Kirkpatrick, Gelatt~Jr, and
  Vecchi]{Kirkpatrick1983optimization}
S.~Kirkpatrick, C.~D. Gelatt~Jr, and M.~P. Vecchi.
\newblock Optimization by simulated annealing.
\newblock \emph{Science}, 220\penalty0 (4598):\penalty0 671--680, 1983.

\bibitem[Koopman(1931)]{Koopman1931}
B.~O. Koopman.
\newblock {H}amiltonian systems and transformation in {H}ilbert space.
\newblock \emph{Proc.~of the National Academy of Sciences U.S.A}, 17\penalty0
  (5):\penalty0 315--318, 1931.

\bibitem[Korda and Mezi{\'c}(2018)]{KordaMezic2018}
M.~Korda and I.~Mezi{\'c}.
\newblock Linear predictors for nonlinear dynamical systems: {K}oopman operator
  meets model predictive control.
\newblock \emph{Automatica}, 93:\penalty0 149--160, 2018.

\bibitem[Liu et~al.(2023)Liu, Ozay, and Sontag]{liu2023properties}
Z.~Liu, N.~Ozay, and E.~D. Sontag.
\newblock Properties of immersions for systems with multiple limit sets with
  implications to learning {K}oopman embeddings, 2023.
\newblock arXiv:2312.17045.

\bibitem[Mauroy et~al.(2020)Mauroy, Susuki, and
  Mezi{\'c}]{MauroySusukiMezic2020}
A.~Mauroy, Y.~Susuki, and I.~Mezi{\'c}.
\newblock \emph{{K}oopman operator in systems and control}.
\newblock Springer, 2020.

\bibitem[Mezi{\'c}(2005)]{Mezic2005}
I.~Mezi{\'c}.
\newblock Spectral properties of dynamical systems, model reduction and
  decompositions.
\newblock \emph{Nonlinear Dynamics}, 41:\penalty0 309--325, 2005.

\bibitem[Mohseni et~al.(2022)Mohseni, McMahon, and Byrnes]{mohseni2022ising}
N.~Mohseni, P.~L. McMahon, and T.~Byrnes.
\newblock {I}sing machines as hardware solvers of combinatorial optimization
  problems.
\newblock \emph{Nature Reviews Physics}, 4\penalty0 (6):\penalty0 363--379,
  2022.

\bibitem[N{\"u}ske et~al.(2023)N{\"u}ske, Peitz, Philipp, Schaller, and
  Worthmann]{NuskePeitzPhilippSchallerWorthmann2023}
F.~N{\"u}ske, S.~Peitz, F.~Philipp, M.~Schaller, and K.~Worthmann.
\newblock Finite-data error bounds for {K}oopman-based prediction and control.
\newblock \emph{Journal of Nonlinear Science}, 33\penalty0 (1):\penalty0 14,
  2023.

\bibitem[Otto and Rowley(2021)]{OttoRowley2021}
S.~E. Otto and C.~W. Rowley.
\newblock {K}oopman operators for estimation and control of dynamical systems.
\newblock \emph{Annual Review of Control, Robotics, and Autonomous Systems},
  4:\penalty0 59--87, 2021.

\bibitem[Proctor et~al.(2018)Proctor, Brunton, and
  Kutz]{ProctorBruntonKutz2018}
J.~Proctor, S.~Brunton, and J.~N. Kutz.
\newblock Generalizing {K}oopman theory to allow for inputs and control.
\newblock \emph{SIAM J.~on Applied Dynamical Syst.}, 17\penalty0 (1):\penalty0
  909--930, 2018.

\bibitem[Rothblum and Tan(1985)]{RothblumTan1985}
U.~G. Rothblum and C.~P. Tan.
\newblock Upper bounds on the maximum modulus of subdominant eigenvalues of
  nonnegative matrices.
\newblock \emph{Linear algebra and its applications}, 66:\penalty0 45--86,
  1985.

\bibitem[Swendsen and Wang(1986)]{Swendsen1986replica}
R.~H. Swendsen and J.-S. Wang.
\newblock Replica {M}onte {C}arlo simulation of spin-glasses.
\newblock \emph{Physical Review Lett.}, 57\penalty0 (21):\penalty0 2607, 1986.

\bibitem[Xu and Antsaklis(2004)]{xu2004optimal}
X.~Xu and P.~Antsaklis.
\newblock Optimal control of switched systems based on parameterization of the
  switching instants.
\newblock \emph{IEEE Trans.~on Automat.~Contr.}, 49\penalty0 (1):\penalty0
  2--16, 2004.

\end{thebibliography}

\end{document}